\documentclass[a4paper,12pt]{article}
\usepackage[english]{babel}
\usepackage[T2A]{fontenc}
\usepackage[cp1251]{inputenc}
\usepackage{amsthm}
\usepackage[tbtags]{amsmath}
\usepackage{amsfonts,amssymb}
\begin{document}

\begin{center}
\textbf{E. I. Kompantseva, A. A. Tuganbaev}

\textbf{Rings on Abelian Torsion-Free Groups of Finite Rank}
\end{center}

\textbf{Abstract.} In the class of reduced Abelian torsion-free groups $G$ of finite rank, we describe $TI$-groups, this means that every associative ring on $G$ is filial. If every associative multiplication on $G$ is the zero multiplication, then $G$ is called a $nil_a$-group. It is proved that a reduced Abelian torsion-free group $G$ of finite rank is a $TI$-group if and only if $G$ is a homogeneous Murley group or $G$ is a $nil_a$-group. We also study the interrelations between the class of homogeneous Murley groups and the class of $nil_a$-groups. For any type $t\ne (\infty,\infty,\ldots)$ и every integer $n>1$, there exist $2^{\aleph_0}$ pairwise non-quasi-isomorphic homogeneous Murley groups of type $t$ and rank $n$ which are $nil_a$-groups. We describe types $t$ such that there exists a homogeneous Murley group of type $t$ which is not a $nil_a$-group.\\
This paper will be published in Beitr\"age zur Algebra und Geometrie / Contributions to Algebra and Geometry.

\textbf{Key words.} Abelian group, ring on Abelian group, filial ring, $TI$-group.

MSC2020 datebase: 20K30, 20K99, 16B99

\section{Introduction}
For an Abelian group $G$, any homomorphism $\mu\colon G\otimes G\to G$ is called a \textsf{multiplication} on $G$. An Abelian group $G$, with multiplication defined on $G$, is called a \textsf{ring on} $G$. On every group $G$, we can define the multiplication $\mu\colon G\otimes G\to 0$ which is called the \textsf{zero multiplication}. A ring with such a multiplication is said to be \textsf{trivial}. If a group $G$ has no  multiplication (resp., no associative multiplication) besides the zero multiplication, then $G$ is called a $nil$-group (resp., a $nil_a$-group). R.A.Beaumont \cite{Be48} posed the description problem for ring structures on Abelian groups; he considered rings on direct sums of cyclic groups.

In the theory of additive groups of rings, one of parts is the study of Abelian groups on which every ring belongs to a certain class. We consider the class of filial rings which were defined in \cite{Er83} and studied in \cite{Sa88},\cite{A-P88}, \cite{F-P04}, \cite{F-P05}. An associative ring is called a \textsf{filial} ring if every meta-ideal of finite index is an ideal. A subring $A$ of an associative ring $R$ is called a \textsf{meta-ideal of index} $n$ if there exists a chain 
$$
A=A_0\subset A_1\subset \cdots\subset A_n=R
$$
such that $A_i$ is an ideal of $A_{i+1}$ for all $i=0,\cdots,n-1$ \cite{Er83}. It is easy to see that an associative ring is filial if and only if the relation `is an ideal of' is transitive. We note that the class of filial rings contains not only Hamiltonian rings, \cite{Red52}, \cite{And67}, \cite{Kru68}, 
but (Von Neumann) regular rings and simple rings; on regular rings, e.g., see\cite{Goo79}, \cite{Tug02}.
In this paper, we study Abelian groups, on which every associative ring is a filial ring; such groups are called $TI$-groups. In \cite{AndW14}, the problem of the study of $TI$-groups is formulated. In the same paper, torsion $TI$-groups are described. In \cite{K-T19}, $TI$-groups are described in the class of algebraically compact Abelian groups. In \cite{KTG21}, $TI$-groups from some classes of Abelian torsion-free groups are described. In \cite{KTG21}, non-reduced $TI$-groups are also described. This allows us to consider only the reduced case in the studies of $TI$-groups.

The main aim of this paper is to describe $TI$-groups in the class of Abelian torsion-free groups of finite rank. In \cite{Mur72} and \cite{Mur74}, C.E.~Murley systematically studied reduced Abelian torsion-free groups $G$ of finite rank such that the $p$-rank of the group $G$ does not exceed $1$, for every prime $p$; these groups are now called \textsf{Murley} groups. We note that L.~Fuchs \cite{Fuc15} decided to break this terminological tradition to pay tribute L.~Prochazka, he calls these groups \textsf{Prochazka-Murley  groups}. Indeed, L.~Prochazka [2] initiated the study of groups of this class a few years before the publication of Murley's paper \cite{Mur72} and obtained more essential results. About Murley groups, see also the book \cite[Section 44]{KMT03}.

Let $\mathcal{MH}$ denote the class of all homogeneous Murley groups and let $\mathcal{NIL}_a$ denote the class of all Abelian torsion-free $nil_a$-groups of finite rank.
 
The main result of Section 2 is that the class of reduced torsion-free groups of finite rank coincides with the class $\mathcal{MH}\,\cup\,\mathcal{NIL}_a$ (Theorem 11). In addition, we prove some properties of Murley groups and rings on them. 

Since the class of reduced torsion-free $TI$-groups of finite rank coincides with the union of the classes $\mathcal{MH}$ and $\mathcal{NIL}_a$, interrelations between these classes are studied in Section 3. It is shown that the classes $\mathcal{MH}\,\cap\,\mathcal{NIL}_a$ and $\mathcal{MH}\,\setminus\,\mathcal{NIL}_a$ are very large and varied.

It is easy to see that every torsion-free group $G$ of rank $1$ belongs to the class $\mathcal{MH}$. In addition, $G\in \mathcal{MH}\,\cap\,\mathcal{NIL}_a$ if the type $t(G)$ is not an idempotent type and $G\in \mathcal{MH}\,\setminus\,\mathcal{NIL}_a$ if the type $t(G)$ is an idempotent type. But the situation changes 
significantly if the rank of the group exceeds $1$. In Theorem 15, it is proved that for any type $t\ne (\infty,\infty,\ldots)$ and each positive integer $n>1$, the class $\mathcal{MH}\,\cap\,\mathcal{NIL}_a$ contains $2^{\aleph_0}$ of pairwise non-quasi-isomorphic groups for type $t$ and rank $n$. 

However, not all types $t\ne (\infty,\infty,\ldots)$ can be realized as types of some groups in the class $\mathcal{MH}\,\setminus\,\mathcal{NIL}_a$. For example, this class does not contain groups of rank $>1$ and type $(0,0,\ldots)$. In Theorem 17, a criterion is obtained for a type to be the type of some group in $\mathcal{MH}\,\setminus\,\mathcal{NIL}_a$. 

In the paper, we widely use a method of the study of rings on torsion-free Abelian groups of finite rank. This argument is based on the method of R.A.Beaumont and R.S.Pierce \cite{B-P61}. For such rings, they proved an analogue of the main Wedderburn theorem on finite-dimensional algebras. They also offered a reasonably simple system of invariants for the quasi-isomorphism classes of quotient divisible torsion-free groups of finite rank.

All groups, considered in the paper, are Abelian, and the word "group"~ means "Abelian group"~ in what follows. 
A multiplication $\mu\colon G\otimes G\to G$ on the group $G$ is often denoted by $\times$, i.e. $ \mu(g_1\otimes g_2)=g_1\times g_2 $ for all $g_1,g_2\in G$. A group $G$ with multiplication $\times$, defined on $G$, determines a ring on the group $G$ which is denoted by $(G,\times)$.
As usual, $\mathbb{N}$, $\mathbb{N}_0$, and $\mathbb{P}$ are the sets of all positive integers, non-negative integers, and prime integers, respectively; $\mathbb{Z}$ is the ring of integers, $\mathbb{Q}$ is the group (or the field) of rational numbers, $\widehat{\mathbb{Z}}_p$ is the group (or the ring) of $p$-adic integers, and ${\mathbb{Q}}_p^*$ is the field of $p$-adic numbers. If $a,b\in \mathbb{Z}$, then $(a,b)$ is the largest common divisor of the integers $a$ and $b$.

Let $G$ be a group, $g\in G$, and let $(G,\times)$ be a ring on $G$. We denote by $(g)_\times$ the ideal of the ring $(G,\times)$ generated by the element $g$. We denote by $\chi(g)$ and $t(g)$ the characteristic and the type of the element $g\in G$, respectively; $t(G)$ is the type of the group $G$ provided $G$ is homogeneous. We denote by $r(G)$ and $r_p(G)$ the rank and the $p$-rank of the group $G$, respectively; $\bar{G}$ is the divisible hull of the group $G$. If $S\subseteq G$, then $\langle S\rangle$ and $\langle S\rangle_*$ are the subgroup and the pure subgroup of the group $G$ generated by the set $S$, respectively.

If $R$ is a ring containing $1$, then $Re$ is the cyclic module over $R$ generated by the element $Re$. We write $A\lhd R$ when $A$ is an ideal of the ring $R$. We denote by $A\oplus B$ and $A\dotplus B$\label{est} the direct sum of groups and the direct sum of rings, respectively. Unless otherwise stated, we use definitions and designations from \cite{Fuc15} and \cite{KMT03}.

\section{Torsion-Free $TI$-Groups of Finite Rank}\label{section2}

The aim of this section is to describe $TI$-groups in the class of reduced torsion-free groups of finite rank.

\textbf{Lemma 1.} Let $n\in\mathbb{N}$,~ $G$ be a reduced torsion-free group, and let $A\oplus B$ be a subgroup of index $n$ of $G$ such that $A$ is not a $nil_a$-group and $B\ne 0$. Then $G$ is not a $TI$-group.

\textbf{Proof.} Let $C=A\oplus B$. By  \cite[Corollary 1]{KTG21}, there exists an associative non-filial ring $(C,\times)$ such that $C\times C\subseteq n^2C$. By \cite[Lemma 5]{KTG21}, $G$ is not a $TI$-group.~$\square$

\textbf{Remark 2.}\label{1not1}
According to \cite[Theorem 1]{A-P88}, the ring $(G,\times)$ is filial if and only if for any $g\in G$, we have
$$
(g)_\times=(g)_\times^2+\mathbb{Z}g.
$$

We note that every trivial ring is filial. In the following theorem, it is proved that there do not exist filial nilpotent rings besides trivial.

\textbf{Theorem 3.} Every filial nilpotent torsion-free ring is trivial.

\textbf{Proof.} Let $(G,\times)$ be a filial ring and let $G^t=0$ for some positive integer $t>2$. We take $x\in G^{t-2}$, $p\in P$, and set $x_1=px\in G^{t-2}$. We have
$$
(x_1)_{\times}=x_1\times G+G\times x_1+\mathbb{Z}x_1.
$$
Since $x_1\times G\times G=G\times x_1\times G=G\times G\times x_1=0,$ we have
$$
(x_1)^2_{\times}+\mathbb{Z}x_1=\mathbb{Z}x_1^2+\mathbb{Z}x_1.
$$
Since the ring $(G,\times)$ is filial, $(x_1)_{\times}=(x_1)^2_{\times}+\mathbb{Z}x_1$ by Remark 2. Consequently,
$$
(x_1)_{\times}=\mathbb{Z}x_1^2+\mathbb{Z}x_1. \eqno (1)
$$
Therefore,
$$
x_1\times x=kx_1^2+mx_1. \eqno (2)
$$
for some $k,m\in\mathbb{Z}$.

Let $g\in G$.

\textsf{Case 1.} Let $m=0$.

It follows from $(2)$ that $px^2=kp^2x^2$, whence $p(kp-1)x^2=0$. Therefore, $x^2=0$, since $G$ is a torsion-free group. Consequently,
$x_1^2=p^2x^2=0.$ 
By $(1)$, we have $(x_1)_{\times}=\mathbb{Z}x_1$. Then
$$
x_1\times g=sx_1 \eqno (3)
$$
for some $s\in \mathbb{Z}$. If $s=0$, then $x_1\times g=0$. 
If $s\ne 0$, then we multiply the both parts of $(3)$ from the right by $g$ and obtain
$$
s(x_1\times g)=x_1\times g\times g=0,\;\text{whence }x_1\times g=0.
$$ 
\textsf{Case 2.} Let $m\ne 0$.

It follows from $(2)$ that $m(x_1\times g)=x_1\times x\times g-kx_1^2\times g=0$. Therefore, $x_1\times g=0$.

Thus, $x\times G=0$ for each $x\in G^{t-2}$. Consequently, $G^{t-1}=G^{t-2}\times G=0$. If $t-1>2$, then $G^{t-2}=0$, by the above. Thus, we have $G^2=0$, i.e., the multiplication $\times$ is the zero multiplication.~$\square$

Let $G$ be a torsion-free group of finite rank, $\bar{G} =\mathbb{Q}\otimes G$ be the divisible hull of the group $G$, and let $(G,\,\times)$ be an associative ring on $G$. Then the multiplication $\times$ is uniquely extended to a multiplication on $\bar{G}$ \cite[Theorem 1.3, p.~675]{Fuc15}. 

By the main Wedderburn theorem, the ring $(\bar{G},\times)$ has a decomposition $\bar{G}=\widetilde{S}\oplus \widetilde{N}$, where $\widetilde{N}$ is a nilpotent ideal of the ring $(\bar{G},\,\times)$ and $\widetilde{S}$ is a semisimple subalgebra of the algebra $(\bar{G},\times)$. By \cite{B-P61}, there exists a subring $S$ and a nilpotent ideal $N$ of the ring $(G,\,\times)$ such that $\overline{S} =\widetilde{S}$, $\bar{N}=\widetilde{N}$, and the subring $S\oplus N$ has a finite index in $(G,\,\times)$.

Let $R$ be an associative algebra over $\mathbb{Q}$. According to \cite{B-P61}, we call a ring $(G,\,\times)$ a \textsf{ring of algebra type $R$} if the induced ring $(\bar{G},\,\times)$ is isomorphic to $R$. In this case, one says that the group $G$ \textsf{admits multiplication of algebra type $R$}.

\textbf{Lemma 4.} Every associative ring on a torsion-free group of finite rank either is nilpotent or contains a subring $C$ of finite index that is of the form $C=M_1\oplus \cdots\oplus M_k\oplus N$, where $N$ is a nilpotent ideal, $k\in\mathbb{N},$ $M_i$ is the full matrix ring of order $r_i$ over a ring of rational division algebra type (in particular, $M_i$ is a ring of rational division algebra type for $r_i=1$).

\textbf{Proof.} Let $G$ be a torsion-free group of finite rank and let $(G,\times )$ be an associative ring. Then the induced ring $(\bar{G},\times)$ also is associative. By the main Wedderburn theorem, $\bar{G}=\widetilde{S}\oplus \widetilde{N}$, where $\widetilde{S}$ is a semisimple algebra and $\widetilde{N}$ is a nilpotent ideal. If $\widetilde{S}=0$, then the ring $(\bar{G},\times)$, whence $(G,\times)$ is nilpotent.

Let $\widetilde{S}\neq 0$. By \cite[Theorem 1.4]{B-P61}, the ring $(G,\times)$ has a subring $S\oplus N$ of finite rank such that $N$ is a nilpotent ideal and $S$ is a ring of semisimple algebra type $\widetilde{S}$. By \cite[Corollary 3.5]{B-P61}, the  ring $S$ has a subring $A_1\dotplus \cdots \dotplus A_2$ of finite index, where $A_i$ is a ring of simple algebra type for every $i=1,\;\ldots,\;k$. By \cite[Corollary 3.9]{B-P61}, for every $i=1,\;\ldots,\;k$, the ring $A_i$ contains a subring $M_i$ of finite index which is a full matrix ring of order $r_i\in\mathbb{N}$ over some ring of rational division algebra type.
	
Thus, the subring $C=M_1\oplus \cdots \oplus M_k\oplus N$ is of finite index in the ring $(G,\times)$.~$\square$

Two Abelian groups $A$ and $B$ are said to be \textsf{quasi-isomorphic} (we write $A\sim B$ in this case) if there exist subgroups $A'\subseteq A$ and $B'\subseteq B$ such that $A'\cong B'$, $mA\subseteq A'$,  and $nB\subseteq B'$ for some positive integers $m$ and $n$. A group $G$ is said to be \textsf{quasi-decomposable} if $G\sim A\oplus B$ for some non-zero groups $A$ and $B$. Groups, that are not quasi-decomposable, are called a \textsf{strongly indecomposable}.

\textbf{Corollary 5.} Let $G$ be a strongly indecomposable torsion-free group of finite rank. Then every associative ring $(G,\times)$ on $G$ is either nilpotent or a ring of rational division algebra type.

\textbf{Proof.} Let us assume that the ring $(G,\times)$ is not nilpotent. By Lemma 4, this ring contains a subring $C$ of finite index such that $C=M_1\oplus \cdots \oplus M_k\oplus N$, where $N$ is a nilpotent ideal, $k\in\mathbb{N},$ $M_i$ is the full matrix ring of order $r_i$ over a ring $R_i$ of rational division algebra type. Then
algebra types of the rings $(G,\times)$ and $C$ coincide. The additive group $C^+$ of the ring $C$ is quasi-isomorphic to the group $G$ and it can be represented in the form $C^+=R_1^+\oplus B$, where $R_1^+$ is the additive group of the ring $R_1$, which is a ring of rational division algebra type and $B$ is a group. If either $N\ne 0$, or $k>1$, or $r>1$, then $B\ne 0$; this contradicts to the property that $G$ is strongly indecomposable. Therefore, the ring $C=R_1$ is a ring of rational division algebra type.  Therefore, the ring $(G,\times)$ is a ring of rational division algebra type, as well.

\textbf{Corollary 6.} Every associative non-trivial ring on $TI$-group is a ring of rational division algebra type.

\textbf{Proof.} Let $G$ be a $TI$-group and let $(G,\times)$ be an associative non-trivial ring. By Theorem 3, this ring is not nilpotent. By Lemma 4, this ring has a subring of finite index $C=M_1\oplus \cdots\oplus M_k\oplus N$, where $k\in\mathbb{N}$ and for every $i=1,\;\ldots,\;k$, the additive group $M_i$ is a direct sum of groups on each of which there exists a ring of rational division algebra type. By Lemma 1, the ring $(G,\times)$  is a ring of rational division algebra type.~$\square$

Further, we prove that the class of reduced torsion-free $TI$-groups of finite rank coincides with the union of the class of Murley groups and the class of $nil_a$-groups. 

We recall that a \textsf{Murley group} is a reduced torsion-free group of finite rank such that $r_p(G)\le 1$ for any $p\in P$. We prove some properties of Murley groups and rings on them.

\textbf{Remark 7.} By \cite{Mur72}, every indecomposable Murley group is strongly indecomposable.

\textbf{Theorem 8.}\label{l7}\\
\textbf{1.} On a Murley group, every multiplication is associative and commutative.

\textbf{2.} On a Murley group, every non-trivial ring is not  nilpotent.

\textbf{3.} Every homogeneous Murley group is strongly indecomposable.

\textbf{Proof.} Let $G$ be a Murley group and let $(G,\times)$ be a ring on $G$. The multiplication $\times$ is uniquely extended to the multiplication on the pure-injective hull $\widehat{G}$ of the group $G$ \cite{Fuc15}. In addition, $\widehat{G}=\prod\limits_{p\in P_1}\widehat{G}_p=\prod\limits_{p\in P_1}\widehat{\mathbb{Z}}_pe_p$, where $P_1\subseteq\mathbb{P}$, $\widehat{G}_p$ is the $p$-adic completion of the group $G$, $e_p\in\widehat{G}_p$. According to \cite{Kom10}, the ring $(\widehat{G},\times)$ is the direct product  of ideals $\widehat{\mathbb{Z}}_pe_p$, $p\in P_1$. In addition, if $e_p\times e_p=\tau_p e_p$, where $\tau_p\in \widehat{\mathbb{Z}}_p,$ then $\alpha e_p\times \beta e_p=(\alpha \beta \tau_p) e_p$ for all $p\in P_1$ and $\alpha, \beta\in \widehat{\mathbb{Z}}_p.$

\textbf{1.} It is easy to see that the ring $(\widehat{G},\times)$ is associative and commutative. Therefore, the subring $(G,\times)$ of $(\widehat{G},\times)$ is associative and commutative.
	
\textbf{2.} Let the ring $(G,\times)$ be non-trivial. Then $a\times b\neq 0$ for some $a=(\alpha_p e_p)_{p\in P_1}$ and $b=(\beta_p e_p)_{p\in P_1}\in G\subseteq\widehat{G},$ where $\alpha_p,\beta_p\in \widehat{\mathbb{Z}}_p$, $p\in P_1$. In this case, $\alpha_q\neq 0$, $\beta_q\neq 0$ and $e_q\times e_q=\tau_q e_q\neq 0$ for some $q\in P_1$. Consequently, $\alpha_q^n\tau_q^{n-1}e_q\neq 0$ for every $n\in \mathbb{N}.$ Since $a^n=(\alpha_p^n\tau_p^{n-1}e_p)_{p\in P_1}$, we have that $a^n\neq 0$ for all $n\in\mathbb{N}.$

\textbf{3.} Let $G$ be a homogeneous Murley group and $G=A\oplus B$. Then there exists a prime integer $p$ with $r_p(G)=1$. Since 
$$
1=r_p(G)=r_p(A)+r_p(B),
$$
exactly one of the numbers $r_p(A)$ and $r_p(B)$ is equal to $1$. For definiteness, let $r_p(A)=1$ and $r_p(B)=0$. This means that $h_p(a)\neq \infty$ and $h_p(b)=\infty$ for some $a\in A\setminus\{0\}$ and some $b\in B$. Since the group $G$ is homogeneous, $B=0$. Consequently, $G$ is strongly indecomposable by Remark 7.~$\square$

 We note that the converse to Theorem 8(3) is not true.

\textbf{Proposition 9.} There exist strongly indecomposable non-homogeneous Murley groups.

\textbf{Proof.} In the group $V=\mathbb{Q}e_1\oplus\mathbb{Q}e_2$, the linear combinations of the form $c=k_1e_1+k_2e_2$, where $k_1,k_2\in\mathbb{Z}$, $k_1>0$, $k_2\ne 0$ and $(k_1,k_2)=1$, will be ordered in the form of a sequence $\{c_j\}_{j>2}$. We set $c_1=e_1$, $c_2=e_2$. Let $\{p_j\}_{j\in \mathbb{N}}$ be the sequence of all prime integers. In $V$, we consider the subgroup $G=\langle p_j^{-\infty}c_j\,|\,j\in\mathbb{N} \rangle$. Then the rank $G$ is equal to $2$. In addition, $G$ is a Murley group. Indeed, let $a$ and $b$ be two linearly independent elements of $G$, $j\in\mathbb{N}$. If $h_{p_j}(a)=\infty$ or $h_{p_j}(b)=\infty$, then $a$, $b$ are $p_j$-dependent. We assume that $h_{p_j}(a)\ne\infty$ and $h_{p_j}(b)\ne\infty$. Then $a$, $b$, $c_j$ are linearly dependent, i.e., $m_1c_j=m_2a+m_3b$ for some $m_1,m_2,m_3$ which are not equal to zero simultaneously. It is easy to see that $m_2\ne 0$ or $m_3\ne 0$. For definiteness, let $m_2\ne 0$. Then $p_j^t\,\not\vert\,m_2$ for some positive integer $t$. However, $p_j^t\,\vert\,m_1c_j$. Consequently, $r_{p_j}(G)\le 1$.\\
In the group $G$, any two linearly independent elements have incomparable types. According to \cite{DeG57}, the group $G$ is pure-indecomposable in this case, i.e., every pure subgroup of the group $G$ is indecomposable. By Remark 7, the group $G$ is strongly indecomposable.~$\square$

In the class of Murley groups, it follows from Theorem 8(1) that the notion of a $nil$-group coincides with the notion of a $nil_a$-group.

\textbf{Theorem 10.} Let $G$ be a Murley group which is not a $nil$-group. Then the following conditions are equivalent.

\textbf{1)} $G$ is a homogeneous group of idempotent type.

\textbf{2)} $G$ is an indecomposable group.

\textbf{3)} $G$ is a strongly indecomposable group.

\textbf{4)} Every non-trivial ring on $G$ is a ring of field type.

\textbf{5)} On the group $G$, there exists a ring of field type.

\textbf{Proof.}
The implication 1)\,$\Rightarrow$\,2) follows from Theorem 8(3).

The implication 2)\,$\Rightarrow$\,3) follows from Remark 7.

The implication 3)\,$\Rightarrow$\,4) follows from Corollary 6 and Theorem 8.

The implication 4)\,$\Rightarrow$\,5) follows from the property that $G$ is not a $nil$-group, by assumption.

5)\,$\Rightarrow$\,1). By \cite[Lemma 9.7]{B-P61}, every group $G$, which admits multiplication of simple algebra type, is homogeneous. The type of $G$ is an idempotent type, since otherwise $G$ is a $nil$-group.~$\square$

We note that the assertion of Theorem 10 is not true without the condition that there exists a non-trivial ring on the group $G$. Indeed, the Murley group $G$, which is constructed in Proposition 9, satisfies condition 3) but does not satisfy condition 1), since $G$ is a $nil$-group.

In the following theorem, we describe $TI$-groups in the class of reduced torsion-free groups of finite rank.

\textbf{Theorem 11.} A reduced torsion-free group $G$ of finite rank is a $TI$-group if and only if $G$ is a homogeneous Murley group or $G$ is a $nil_a$-group

\textbf{Proof.} Let $G$ be a $TI$-group such that there exists a non-trivial associative ring $(G,\times)$. Then $(G,\times)$ is a ring of rational division algebra type, by Corollary 6. By \cite[Lemma 9.7]{B-P61}, the group $G$ is homogeneous.

We assume that $r_p(G)\ge 2$ for some $p\in P$. Then $G$ contains $p$-independent elements $a$ and $b$. Since the induced ring $(\bar{G},\times)$ on the divisible hull $\bar{G}$ of the group $G$ is a rational division algebra, $(\bar{G},\times)$ contains the inverse element $a^{-1}$ of $a$. In addition, $na^{-1}\in G$ for some positive integer $n$. Let $k$ be a positive integer such that $p^k$ does not divide $n$. In the ring $(G,\times)$, the ideal $(p^ka)_{\times}$ satisfies the relation
$$
(p^ka)_{\times}^2+\mathbb{Z}p^ka\subseteq p^{2k}G+\mathbb{Z}p^ka.
$$
We consider the element $p^knb=p^ka\times na^{-1}\times b\in (p^ka)_{\times}$. We assume that $p^knb\in (p^ka)_{\times}^2+ \mathbb{Z}p^ka$. Then $p^knb=p^{2k}x+sp^ka$ for some $x\in G$, $s\in\mathbb{Z}$. Consequently, $nb-sa=p^kx$. It follows from $p$-independence of the elements $a$ and $b$ that $p^k$ divides $n$; this contradicts to the choice of the integer $k$. Therefore, $(p^ka)_{\times}\ne (p^ka)_{\times}^2+\mathbb{Z}p^ka$. Therefore, the ring $(G,\times)$ is not filial. Consequently, $r_p(G)\le 1$ for every $p\in P$.

It is easy to see that every $nil_a$-group is a $TI$-group. 
Let $G$ be a homogeneous Murley group and there exists a non-trivial ring $(G,\times)$ on $G$. Then $t(G)$ is an idempotent type. By Theorem 10, the ring $(G,\times)$ is a ring of field type. In this case, it follows from \cite[Theorem 7]{AndS07} that the ring $(G,\times)$ is filial if for any $g\in G$, there exists a positive integer $m$ such that $mg\in G\times g^2$. Let $g\in G$. Since the induced ring $(\bar{G},\times)$ is a field, it contains the inverse element $g^{-1}$ of $g$. In addition, there exists a positive integer $m$ such that $mg^{-1}\in G$, whence 
$mg=mg^{-1}\times g^2\in G\times g^2.$ Consequently, the ring $(G,\times)$ is filial. Since the ring $(G,\times)$ is arbitrary, $G$ is a $TI$-group.~$\square$

\textbf{Corollary 12.} A reduced torsion-free group $G$ of finite rank is a $TI$-group if and only if either $G$ is a $nil_a$-group or $G$ is a Murley group which satisfies one of conditions \textbf{1)}--\textbf{5)} of Theorem 10.

\section{Homogeneous Murley Groups}\label{section3}

Let $\mathcal{MH}$ be the class of all homogeneous Murley groups and let $\mathcal{NIL}$ (resp., $\mathcal{NIL}_a$) be the class of all torsion-free $nil$-groups (resp., $nil_a$-groups) of finite rank.

According to Theorem 11, the class of all torsion-free $TI$-groups of finite rank coincides with $\mathcal{MH}\,\cup\, \mathcal{NIL}_a$. However, this theorem does not explain the interrelations between the classes $\mathcal{MH}$ and $\mathcal{NIL}_a$. For example, what can we say about classes 
$$
\mathcal{MH}\,\cap\,\mathcal{NIL}_a \qquad\text{and}\qquad \mathcal{MH}\setminus \mathcal{NIL}_a\;?
$$
In particular, what types and positive integers can be realized as types and ranks of groups in these classes provided the classes are not empty?

We note that
$$
\mathcal{MH}\,\cap\,\mathcal{NIL}_a=\mathcal{MH}\,\cap\,\mathcal{NIL}, \quad \mathcal{MH}\setminus \mathcal{NIL}_a=\mathcal{MH}\setminus \mathcal{NIL}
$$
by Theorem 8(1).

If $r(G)=1$, then $G\in \mathcal{MH}$; in addition $G\in \mathcal{MH}\,\cap\,\mathcal{NIL}$ if and only if the type $t(G)$ is not an idempotent type.

It follows from Theorem 10 that  there exists a ring of field type on every group in the class $\mathcal{MH}\setminus \mathcal{NIL}$. In \cite{B-P61}, it is proved that every group, which admits a multiplication of simple algebra type, is quotient divisible. According to \cite{B-P61}, a torsion-free group $G$ of finite rank is said to be \textsf{quotient divisible} if $G$ has a free Abelian subgroup $F$ such that $G/F$ is a divisible torsion group.
Now quotient divisible groups are actively studied by various authors; e.g., see \cite{AlbBV07}, \cite{Fom09}, \cite{KomF19} and others.

In \cite{B-P61}, there is a quite simple independent invariant system for quasi-isomorphism classes of quotient divisible groups. These invariants are constructed similar to Kurosh, Malcev and Derry invariants for isomorphism classes of arbitrary torsion-free groups of finite rank; they are based on the following argument.

Let $G$ be a torsion-free group of finite rank.
We set $V=\mathbb{Q}\otimes G=\bar{G}$, $V^{(p)}={\mathbb{Q}}_p^*\otimes G$. We assume that the group $V$ is embedded in the group $V^{(p)}$. If $R$ is a subgroup of the additive group ${\mathbb{Q}}_p^*$ containing $\mathbb{Z}$ and $A$ is a subgroup of the group $V^{(p)}$, then we follow \cite{B-P61} and \cite{Mur72} and, instead of $R\otimes A$, we write
$$
RA=\{r_{1}a_{1}+\ldots+r_{k}a_{k}\,|\,r_{i}\in R,\,a_{i}\in A,\, i=1,\;\ldots,\;k\}.
$$ 
If $S\subseteq V^{(p)}$ and $A=\langle S\rangle_R$, then the group $RA$ is also denoted by $\langle S\rangle_R$.

Let $\delta_p(G)$ be the maximal divisible subgroup of the group $\widehat{\mathbb{Z}}_pG$. This group is the maximal divisible submodule of the $\widehat{\mathbb{Z}}_p$-module $\widehat{\mathbb{Z}}_pG$. Therefore, $\delta_p(G)$ is a ${\mathbb{Q}}_p^*$-subspace of the space $V^{(p)}$.

For every $p\in P$, we denote by $\mathcal{L}_p$ the subspace lattice of the ${\mathbb{Q}}_p^*$-space $V^{(p)}$. Let $\mathcal{L}=\prod\limits_p\mathcal{L}_p$ be the direct product of these lattices. Then $\delta_p(G)\in \mathcal{L}_p$ for every $p\in P$. Following \cite{B-P61}, we call $\delta(G)=(\delta_p(G))_{p\in P}\in\mathcal{L}$ the \textsf{quotient divisible invariant} (\textsf{q.d.-invariant}) of the group $G$.

Below we use the following definition everywhere. 

Let $P_1\subseteq P$ and $n\in\mathbb{N}$. We say that the set $\{\alpha_p\in\widehat{\mathbb{Z}}_p\,\vert\,p\in P_1\}$ \textsf{satisfies condition $(T)$ for $n\in \mathbb{N}$} if for any polynomial $f(x)\in \mathbb{Z}[x]$ of degree $\le n$, the following two conditions are true.

\textbf{1)} $f(\alpha_p)\ne 0$ for all $p\in P_1$;

\textbf{2)} $f(\alpha_p)\not\equiv 0\text{ (mod }p)$ for almost all $p\in P_1$.

For the set $\{\alpha_p\in\widehat{\mathbb{Z}}_p\,\vert\,p\in P_1\}$, the elements $\alpha_p$, $p\in P_1\;$, are called \textsf{$p$-components} of this set. It is easy to see that if the set $\{\alpha_p\in\widehat{\mathbb{Z}}_p\,\vert\,p\in P_1\}$ satisfies condition $(T)$ for $n>1$, then it also satisfies condition $(T)$ for any positive integer $m<n$.

If $t$ is a type, then we set
$$
P_{\infty}(t)=\{p\in P\,\vert\,t(p)=\infty\},\;
P_{0}(t)=\{p\in P\,\vert\,t(p)<\infty\}.
$$

\textbf{Lemma 13.} Let $n\in\mathbb{N}$, $n>1$, $t$ be a type which is not equal to $(\infty,\infty,\ldots)$, and let the set $\{\alpha_p\in\widehat{\mathbb{Z}}_p\,\vert\,p\in P_0(t)\}$ satisfy condition $(T)$ for $n-1$. Then there exists a group $G\in \mathcal{MH}$ of rank $n$ and type $t$ such that for some basis $e_1,\;\ldots,\; e_n$ of the space $\bar{G}$, the quotient divisible invariant $\delta(G)=(\delta_p(G))_{p\in P}$ has the form 
$$
\delta_p(G)=V^{(p)} \text{ for } p\in P_{\infty}(t),\;
\delta_p(G)=\langle u_{p,1},\;\ldots,\;u_{p,n-1}\rangle_{\mathbb{{Q}_p^{*}}} \text{ for } p\in P_{0}(t),
$$
$$
\text{where }u_{p,1}=e_1+\alpha_pe_2,\; \ldots,\;u_{p,n-1}=e_{n-1}+\alpha_pe_n.
$$
In addition, if the type $t$ is an idempotent type, then the group $G$ is quotient divisible.

\textbf{Proof.} We set $P_{\infty}=P_{\infty}(t)$, $P_{0}=P_{0}(t)$. Let $p\in P_0$. We write the number $\alpha_p$ in the form
$$
\alpha_p=a_{p,0}+a_{p,1}p+a_{p,2}p^2+\ldots, \text{ where } a_{p,k}\in \mathbb{Z},\; 0\le a_{p,k}\le p-1.
$$
For every $k\in\mathbb{N}$, we set
$$
\alpha_{p,k}=a_{p,0}+a_{p,1}p+\ldots+a_{p,k-1}p^{k-1}.
$$
We take a characteristic $\chi=(\chi_p)_{p\in P}\in t$. For all $p\in P_0$, $s\in \{1,\;\ldots,\;n-1\}$ and $k\in\mathbb{N}$, we consider elements
$$
u_{p,s}^{(k)}=e_{s}+\alpha_{p,k}e_{s+1},\qquad 
x_{p,s}^{(k)}=p^{-(k+\chi_p)}u_{p,s}^{(k)} \eqno (4)
$$
of the space $V=\mathbb{Q}e_1\oplus\ldots\oplus\mathbb{Q}e_n$
Let $R$ be a subgroup of the group $\mathbb{Q}$ such that $1\in R$ and $\chi_R(1)=\chi$. In the space $V$, we consider the subgroups
$$
F=\mathbb{Z}e_1\oplus\ldots\oplus \mathbb{Z}e_n,\; 
A=Re_1\oplus\ldots\oplus Re_n,\;
$$
$$
G=\langle\,\{x_{p,s}^{(k)}\,\vert\,p\in P_0,\, s\in \{1,\;\ldots,\;n-1\},\, k\in\mathbb{N}\},\, A\rangle_R,
$$
$$
A_{p,k}=\langle x_{p,1}^{(k)},\;\ldots,\;x_{p,n-1}^{(k)},e_n\rangle_R,\quad p\in P_0,\quad k\in \mathbb{N}.
$$
Let $p\in P_0$, $s\in \{1,\;\ldots,\;n-1\}$, $k\in\mathbb{N}$. By $(4)$, we have $e_{s}=p^{k+\chi_p}x_{p,s}^{(k)}-\alpha_{p,k}e_{s+1}$. Consequently,
$$
A\subseteq \langle x_{p,1}^{(k)},\;\ldots,\;x_{p,n-1}^{(k)},e_n\rangle_{R}= A_{p,k} \text{ for all } p\in P_0,\,k\in\mathbb{N}. \eqno (5)
$$
If $m\in\mathbb{N}$ and $m<k$, then
$$
p^{m+\chi_p}x_{p,s}^{(m)}=e_s+\alpha_{p,m}e_{s+1},\;
p^{k+\chi_p}x_{p,s}^{(k)}=e_s+\alpha_{p,k}e_{s+1}.
$$
Therefore, $p^{m+\chi_p}x_{p,s}^{(m)}-\alpha_{p,m}e_{s+1}=p^{k+\chi_p}x_{p,s}^{(k)}-\alpha_{p,k}e_{s+1},$
$$
p^{m+\chi_p}x_{p,s}^{(m)}=p^{k+\chi_p}x_{p,s}^{(k)}-(\alpha_{p,k}-\alpha_{p,m})e_{s+1}.
$$
Then~ $x_{p,s}^{(m)}=p^{k-m}x_{p,s}^{(k)}-(a_{p,m}+\ldots+a_{p,k-1}p^{k-m-1})p^{-\chi_p}e_{s+1}.$
Consequently,
$$
x_{p,s}^{(m)}\in\langle x_{p,s}^{(k)},\;\ldots,\; x_{p,n-1}^{(k)}, e_n\rangle_R =A_{p,k}\, \text{for all } p\in P_0,\, m<k. \eqno (6)
$$

Now let $q\in P_0$, $\;q\ne p$, $\;k,m\in\mathbb{N}$. Since
$$
q^{m+\chi_q}x_{q,s}^{(m)}=e_s+\alpha_{q,m}e_{s+1}\, \text{and }\,
p^{k+\chi_p}x_{p,s}^{(k)}=e_s+\alpha_{p,k}e_{s+1},
$$
we have that
$$
q^{m+\chi_q}x_{q,s}^{(m)}=p^{k+\chi_p}x_{p,s}^{(k)}-
\alpha_{p,k}e_{s+1}+\alpha_{q,m}e_{s+1}=
p^{k+\chi_p}x_{p,s}^{(k)}+(\alpha_{q,m}-\alpha_{p,k})e_{s+1}.
$$
Therefore, it follows from $(5)$ that 
$$
q^{m+\chi_q}x_{q,s}^{(m)}\in
\langle x_{p,1}^{(k)},\;\ldots,\; x_{p,n-1}^{(k)}, e_n\rangle_R =A_{p,k} \text{ for all } p,q\in P_0,\, k,m\in \mathbb{N}.
\eqno (7)
$$
Let $g\in G$. By $(5)$ and $(6)$, the element $g$ can be represented in the form $g=\sum_{q\in P_1}g_q$, where $P_1$ is some finite subset in $P_0$ and $g_q\in A_{q,m_q}$ for some $m_q\in\mathbb{N}$.  We fix $p\in P_0$. It follows from $(7)$ that there exist integer $b,k\in\mathbb{N}$ such that
$$
bg\in A_{p,k}\, \text{ and }\, (b,p)=1. \eqno (8)
$$
We prove that the subgroup $Re_n$ is pure in $G$. Let $re_n=p^mg$ for some $p\in P_0$, $m\in \mathbb{N}$, $r\in R$, and $g\in G$. Then the element $g$ satisfies condition $(8)$ for some $b\in\mathbb{N}$, $k\in\mathbb{N}$. Therefore, 
$$
bre_n=bp^mg=p^m(c_1x_{p,1}^{(k)}+\ldots + c_{n-1}x_{p,n-1}^{(k)}+c_{n}e_n),
$$
where $c_{1},\ldots ,c_{n}\in R$.
By substituting  $x_{p,1}^{(k)},\;\ldots,\;x_{p,n-1}^{(k)}$ from $(4)$ to the right part of this relation, we obtain that the coefficient of $e_1$ is equal to $p^mc_1$. Therefore, $c_1=0$. Similarly, $c_2=\ldots c_{n-1}=0$. Thus, $bre_n=p^mc_ne_n$, whence $re_n=p^mc_n'e_n$, where $c_n'\in R$. Consequently, the subgroup $Re_n$ is $p$-pure for every $p\in P_0$. Since $Re_n$ is $p$-divisible for all $p\in P_{\infty}$, we have that $Re_n$ is a pure subgroup of the group $G$.

We prove that $G$ is a homogeneous group of type $t$. Let $a=r_{1}e_{1}+\ldots+r_{n}e_{n}\in F$, where $r_i\in \mathbb{Z}$. It is easy to see that
$$
\chi (a)\ge \chi, \eqno (9)
$$
since $\chi (e_i)=\chi$ for all $i\in \{1,\;\ldots,\;n\}$. In particular, $h_p(a)=\infty$ for all $p\in P_{\infty}$.

Let $p\in P_{0}$, $k\in\mathbb{N}$. In the relations below, the symbol $\equiv$ means congruence modulo $p^{k+\chi_p}$. First, we note that
$$
a=r_1e_1+\ldots+r_ne_n\equiv r_1e_1+\ldots+r_ne_n-r_1p^{k+\chi_p}x_{p,1}^{(k)}=
$$
$$
=r_1e_1+\ldots+r_ne_n-r_1(e_1+\alpha_{p,k}e_2)= (r_2-r_1\alpha_{p,k})e_2+r_3e_3\ldots+r_ne_n,
$$
and for any $s\in \{2,\;\ldots,\;n-1\}$, it follows from $(4)$ that $$
(r_s-r_{s-1}\alpha_{p,k}+r_{s-2}\alpha_{p,k}^2+\ldots+(-1)^{s-1}r_1\alpha_{p,k}^{s-1})e_s+r_{s+1}e_{s+1}+\ldots+r_ne_n\equiv
$$
$$
\equiv (r_s-r_{s-1}\alpha_{p,k}+r_{s-2}\alpha_{p,k}^2+\ldots+(-1)^{s-1}r_1\alpha_{p,k}^{s-1})e_s+r_{s+1}e_{s+1}+\ldots+r_ne_n-
$$
$$
-(r_s-r_{s-1}\alpha_{p,k}+r_{s-2}\alpha_{p,k}^2+\ldots+
(-1)^{s-1}r_1\alpha_{p,k}^{s-1})p^{k+\chi_p}x_{p,s}^{(k)}=
$$
$$
r_{s+1}e_{s+1}+\ldots+r_ne_n-
(r_s-r_{s-1}\alpha_{p,k}+\ldots+(-1)^{s-1}r_1\alpha_{p,k}^{s-1})\alpha_{p,k}e_{s+1}=
$$
$$
=(r_{s+1}-r_s\alpha_{p,k}+r_{s-2}\alpha_{p,k}^2+\ldots+
(-1)^{s}r_1\alpha_{p,k}^{s})e_{s+1}+r_{s+2}e_{s+2}+\ldots+ r_ne_n.
$$
Consequently, $a=r_1e_1+\ldots+r_ne_n\equiv $
$$
\equiv (r_n-r_{n-1}\alpha_{p,k}+r_{n-2}\alpha_{p,k}^2+\ldots+(-1)^{n-1}r_1\alpha_{p,k}^{n-1})e_n
\, \text{ (mod }p^{k+\chi_p}).\eqno (10)
$$
We set $f(x)=r_n-r_{n-1}x+r_{n-2}x^2+\ldots+(-1)^{n-1}r_1x^{n-1}\in\mathbb{Z}[x].$
Then the relation $(10)$ can be given in the form
$$
a\equiv f(\alpha_{p,k})e_n\;(\text{mod }p^{k+\chi_p}). \eqno (11)
$$
We assume that $p^k$ divides $a$ for every $k\in\mathbb{N}$. Then $f(\alpha_{p,k})e_n\equiv 0 \text{ (mod }p^{k+\chi_p})$ for every $k\in\mathbb{N}$ by $(11)$. Since the subgroup $Re_n$ is pure in $G$, we have that $f(\alpha_{p,k})\equiv 0 \text{ (mod }p^{k})$ for every $k\in\mathbb{N}$. Therefore, the sequence of integers $\{f(\alpha_{p,k})\}_{k\in\mathbb{N}}$ is a Cauchy sequence in the $p$-adic topology on $\widehat{\mathbb{Z}}_p$ and the limit $\lim\{f(\alpha_{p,k})\}_{k\in\mathbb{N}}$ of this sequence is equal to zero. Consequently,
$$
0=\lim\{f(\alpha_{p,k})\}_{k\in\mathbb{N}}=
f(\lim\{\alpha_{p,k}\}_{k\in\mathbb{N}})=f(\alpha_p).
$$
This contradicts to the property that the set $\{\alpha_{p}\,\vert\,p\in P_0\}$ satisfies condition $(T)$ for $n-1$. Therefore, $h_p(a)<\infty$ for all $p\in P_0$.

Let $P_2=\{p\in P_0\,\vert\,p^{\chi_p+1}|a\}$. If $p\in P_2$, then
$$
0\equiv a\equiv f(\alpha_{p,\chi_p+1})e_n \text{ (mod }p^{\chi_p+1})
$$
by $(11)$. Therefore, $f(\alpha_{p,\chi_p+1})\equiv 0\text{ (mod }p)$ for all $p\in P_2$. Since the set $\{\alpha_{p}\,\vert\,p\in P_0\}$ satisfies condition $(T)$ for $n-1$, we have that $P_2$ is a finite set. Therefore, $h_p(a)=\chi_p$ for almost all $p\in P_0$. Therefore, $t(a)=t$ by $(9)$. Since $G/F$ is a torsion group, $G$ also is a homogeneous group of type $t$.

For every $p\in P$, we consider the divisible part $\delta_p(G)$ of the group $\widehat{\mathbb{Z}}_pG$. If $p\in P_{\infty}$, then the group $\widehat{\mathbb{Z}}_pG$ is divisible, whence
$$
\delta_p(G)=\mathbb{Q}_p^*e_1\oplus\ldots\oplus\mathbb{Q}_p^*e_n=V^{(p)}.
$$
Let $p\in P_0$. Since the group $\widehat{\mathbb{Z}}_pG$ is not  divisible, we have \\
$\text{dim}_{\mathbb{Q}_p^*}\delta_p(G)\le n-1$. We consider elements
$$
u_{p,1}=e_1+\alpha_pe_2,\;\ldots,\;u_{p,n-1}=e_{n-1}+\alpha_pe_n
$$
of the group $\widehat{\mathbb{Z}}_pG$. It is easy to see that $u_{p,1},\;\ldots,\;u_{p,n-1}$ are linearly independent over $\mathbb{Q}_p^*$. In addition, the pure subgroup $\langle u_{p,1},\;\ldots,\;u_{p,n-1}\rangle_*$ of the group $\widehat{\mathbb{Z}}_pG$ is divisible, since
$$
u_{p,s}=e_{s}+\alpha_pe_{s+1}=e_{s}+\alpha_{p,k}e_{s+1}+
(\alpha_{p}-\alpha_{p,k})e_{s+1}=
$$
$$
=p^{k+\chi_p}x_{p,s}^{(k)}+(a_{p,k}p^{k}+a_{p,k+1}p^{k+1}+\ldots)e_{s+1}=
$$
$$
=p^{k}(p^{\chi_p}x_{p,s}^{(k)}+(a_{p,k}+a_{p,k+1}p+\ldots)e_{s+1})
$$
for all $s\in\{1,\;\ldots,\;n-1\}$ and $k\in\mathbb{N}$. Consequently, $\delta_p(G)=\langle u_{p,1},\,\;\ldots,\;\,u_{p,n-1}\rangle_{Q_p^*}$.

According to \cite{Mur72}, the relation $r_p(G)+\text{dim }_{\mathbb{Q}_p^*}\delta_p(G)=n$ is true for every $p\in P$. Therefore, $r_p(G)=0$ for $p\in P_{\infty}$ and $r_p(G)=1$ for $p\in P_0$. Consequently, $G$ is a Murley group.

It is easy to see that 
$$
G/A=\langle x_{p,s}^{(k)}+A\;|\;p\in P_0,\,s\in\{1,\;\ldots,\;n-1\},\,k\in\mathbb{N}\rangle
$$
is a divisible torsion group. If $t$ is an idempotent type, then $R/\mathbb{Z}$ is a divisible torsion group; this implies that $A/F$ is a divisible torsion group. Since $G/A\cong (G/F)/(A/F)$, we have that $G/F$ also is a divisible torsion group. Consequently, the group $G$ is quotient divisible.~$\square$

\textbf{Lemma 14.} For every $P_1\subseteq P$ and each $n\in \mathbb{N}$, there exist $2^{\aleph_0}$ sets $\{\alpha_p\in\widehat{\mathbb{Z}}_p\,|\,p\in P_1\}$ such that the following properties are true.\\
The sets $\{\alpha_p\in\widehat{\mathbb{Z}}_p\,|\,p\in P_1\}$  satisfy condition $(T)$ for $n$.\\
Any two of these sets have algebraically independent $p$-components for at least one $p\in P_1$.

\textbf{Proof.} Let $\{p_k\}_{k\in \mathbb{N}}$ be the sequence of all prime integers ordered in ascending order. Let $B_n$ be the set of all polynomials in $\mathbb{Z}[x]$, whose degrees do not exceed $n$ and the coefficients are coprime in total, and let $\{f_k(x)\}_{k\in \mathbb{N}}$ be the sequence of all polynomials in $B_n$ which is ordered so that the leading coefficient of the polynomial $f_k(x)$ is less than $p_k$, $k\in \mathbb{N}$.

By \cite{FinR07}, there exists an integer $c\in \mathbb{N}$ such that $p_k > c k \ln k$ for all $k\in \mathbb{N}$. It is easy to see that $c \ln k_0> n+1$ for some $k_0\in \mathbb{N}$.

Let $k\ge k_0$. The congruence
$$
f_1(x)\cdot\ldots\cdot f_k(x)\equiv 0 \text{ (mod }p_k) \eqno (12)
$$
has at most $\text{deg }(f_1(x)\cdot\ldots\cdot f_k(x))$ solutions provided the leading coefficient of the polynomial $f_1(x)\cdot\ldots\cdot f_k(x)$ is not divided by $p_k$. If $m_s$ is the leading coefficient polynomial $f_s(x)$, $s\in \{1,\;\ldots,\;k\}$, then $m_s<p_k$. Therefore, the leading coefficient of the polynomial $f_1(x)\cdot\ldots\cdot f_k(x)$, which is equal to $m_1\cdot\ldots\cdot m_k$, is not divided by $p_k$. In addition,
$$
\text{deg }(f_1(x)\cdot\ldots\cdot f_k(x))\le k n< k(c\ln k_0-1)\le k c \ln k -k<p_k-k.
$$
Consequently, there exists a positive integer $a_k\in \{1,\ldots,p_k-1\}$ which is not a solution of comparison $(12)$.

By choosing arbitrary positive integers $a_k\in \{1,\ldots,p_k-1\}$ for $k<k_0$, we obtain a sequence $\{a_k\}_{k\in \mathbb{N}}$ such that $a_k$ is not a solution of $(12)$ for almost all $k\in\mathbb{N}$.
For every $k\in \mathbb{N}$, we take an element $\beta_k\in\widehat{\mathbb{Z}}_{p_k}$ which is transcendental over $\mathbb{Q}$ and we set $\alpha_{p_k}=a_k+p_k\beta_k$. Then the element $\alpha_{p_k}$ is transcendental over $\mathbb{Q}$ for every $k\in \mathbb{N}$.

We prove that the set $\{\alpha_p\,\vert\,p\in P_1\}$ satisfies property \textbf{2)} of condition $(T)$ for $n$. Let 
$f(x)\in \mathbb{Z}[x]$, deg$\,f(x)\le n$. Then $f(x)=f_{k_1}(x)$ for some $k_1\in \mathbb{N}$. Let $k_2=\text{max }\{k_0,k_1\}$. For for all $k\ge k_2$, we have
$$
f_1(\alpha_{p_k})\cdot\ldots\cdot f_k(\alpha_{p_k})\equiv
f_1(a_k)\cdot\ldots\cdot f_k(a_k)\not\equiv 0\, \text{ (mod }p_k)
$$
by the choice of the numbers $a_k$. Consequently, 
$$
f(\alpha_{p_k})=f_{k_1}(\alpha_{p_k})\not\equiv 0 \text{ (mod }p_k).
$$
Therefore, the set $\{\alpha_p\,\vert\,p\in P\}$ satisfies condition $(T)$ for $n$. Therefore, the set $\{\alpha_p\,\vert\,p\in P_1\}$ also satisfies condition $(T)$ for $n$.

Let $N_1=\{k\in\mathbb{N}\;|\;p_k\in P_1\}$. We assume that for every $k\in N_1$, we have chosen transcendental elements $\beta_k'\in \widehat{\mathbb{Z}}_{p_k}$ such that $\beta_m$ and $\beta_m'$ are algebraically independent for some $m\in \mathbb{N}_1$. We set $\alpha_{p_k}'=a_k+p_k\beta_k'$, $k\in \mathbb{N}_1$. By the above, the set $\{\alpha_p'\,\vert\,p\in P_1\}$ satisfies condition $(T)$; in addition $\alpha_{p_m}$ and $\alpha_{p_m}'$ are algebraically independent.

For every $p\in P_1$, the group $\mathbb{Q}_p^*$ has $2^{\aleph_0}$ transcendental algebraically independent elements. Therefore, there exist $2^{\aleph_0}$ sets $\{\alpha_p\in \widehat{\mathbb{Z}}_p\,\vert\,p\in P_1\}$ which satisfy condition $(T)$ for $n$ and each two of them have algebraically independent $p$-components for at least one $p\in P_1$.~$\square$

\textbf{Theorem 15.} For every $n\in \mathbb{N}$, $n>1$, and every type $t\ne (\infty,\infty,\ldots)$, the class $\mathcal{MH}\,\cap\,\mathcal{NIL}$ contains $2^{\aleph_0}$ pairwise non-quasi-isomorphic groups from of rank $n$ and type $t$.

\textbf{Proof.} Let $n\in\mathbb{N}$, $n>1$, $t$ be a type not equal to $(\infty,\infty,\ldots)$, $P_{\infty}=P_{\infty}(t)$, and let $P_{0}=P_{0}(t)$. By Lemma 14, there exists a set $\{\alpha_p\in\widehat{\mathbb{Z}}_p\,\vert\,p\in P_0\}$ which satisfies condition $(T)$ for $n$.

By Lemma 13, there exists a group $G\in\mathcal{MH}$ of rank $n$ and type $t$. In addition, there exists a basis $e_1,\;\ldots,\; e_n$ of the space $V=\bar{G}$ such that
$$
\delta_p(G)=\begin{cases}
\widehat{\mathbb{Z}}_pG=\mathbb{Q}_p^*G \text{ for }p\in P_{\infty}\\
\langle u_{p,1},\;\ldots,\;u_{p,n-1}\rangle_{\mathbb{Q}_p^*} \text{ for }p\in P_{0},
\end{cases}
$$ 
where $u_{p,s}=e_s+\alpha_pe_{s+1}$,~ $s\in\{1,\;\ldots,\;n-1\}$.

Let the set $\{\gamma_p\in \widehat{\mathbb{Z}}_p\,|\, p\in P_{0}\}$ also satisfy condition $(T)$ for $n$. By Lemma 13, there also exists a group $G'\in \mathcal{MH}$ of rank $n$ and type $t$ such that for some basis $e_1',\;\ldots,\; e_n'$ of the space $\overline{G'}$, the spaces $\delta_p(G')$, $p\in P_0$ have the form
$$
\delta_p(G')=\begin{cases}
\widehat{\mathbb{Z}}_pG'=\mathbb{Q}_p^*G' \text{ for }p\in P_{\infty}\\
\langle v_{p,1},\;\ldots,\;v_{p,n-1}\rangle_{\mathbb{Q}_p^*} \text{ for }p\in P_{0},
\end{cases}
$$ 
where $v_{p,s}=e_s'+\gamma_pe_{s+1}'$, $s\in\{1,\;\ldots,\;n-1\}$.

Let $\varphi\colon G\to G'$ be a group homomorphism and 
$$
\varphi(e_s)=c_{s,1}e_1'+\ldots+c_{s,n}e_n',\;
c_{s,i}\in\mathbb{Q},\; s,i\in \{1,\;\ldots,\;n-1\}.
$$
For any $p\in P$, we have a homomorphism 
$$
\varphi^{(p)}=\text{id}_{\mathbb{Q}_p^*}\otimes \varphi\colon \widehat{\mathbb{Z}}_pG\to \widehat{\mathbb{Z}}_pG .
$$
Let $p\in P_{0}$, $s\in \{1,\;\ldots,\;n-1\}$. Since $\varphi^{(p)}(\delta_p(G))\subseteq \delta_p(G')$, we have
$$
\varphi^{(p)}(u_{p,s})=m_1v_{p,1}+\ldots+m_{n-1}v_{p,n-1}=
$$
$$
=m_{1}e_1'+(m_{1}\gamma_p+m_2)e_2'+\ldots+(m_{n-2}\gamma_p+m_{n-1})e_n'+m_{n-1}\gamma_pe_n',
$$
where $m_i\in \mathbb{Q}_p^*$, $i\in \{1,\;\ldots,\;n-1\}$. On the other hand,
$$
\varphi^{(p)}(u_{p,s})=\varphi^{(p)}(e_{s}+\alpha_pe_{s+1})=
\varphi(e_{s})+\alpha_p\varphi(e_{s+1})=
$$
$$
=\sum_{i=1}^{n}c_{s,i}e_i'+ \alpha_p\sum_{i=1}^{n}c_{s+1,i}e_i'=b_1e_1'+\ldots+b_ne_n',
$$
where $b_i=c_{s,i}+\alpha_pc_{s+1,i}$.

By equating the coefficients of $e_1',\;\ldots,\;e_n'$ in these relations, we obtain
$$
m_1=b_1,\; m_1\gamma_p+m_2=b_2,\; \ldots,\; m_{n-2}\gamma_p+m_{n-1}=b_{n-1},\quad m_{n-1}\gamma_p=b_n.\eqno (13)
$$
It follows from $(13)$ that
$$
0=b_n-m_{n-1}\gamma_p=b_n-(b_{n-1}-m_{n-2}\gamma_p)\gamma_p=
b_{n}-b_{n-1}\gamma_p+m_{n-2}\gamma_p^2=
$$
$$
=b_{n}-b_{n-1}\gamma_p+b_{n-2}\gamma_p^2+\ldots+(-1)^{k}b_{n-k}\gamma_{p}^{k}+\ldots+(-1)^{n-1}b_{1}\gamma_{p}^{n-1}=
$$
$$
=(c_{s,n}+\alpha_pc_{s+1,n})-(c_{s,n-1}+\alpha_pc_{s+1,n-1})\gamma_p+\ldots+(-1)^{n-1}(c_{s,1}+\alpha_pc_{s+1,1})\gamma_{p}^{n-1}.
$$
Therefore, $\alpha_p$ and $\gamma_{p}$ satisfy the relation
$$
c_{s,n}+c_{s+1,n}\alpha_p-c_{s,n-1}\gamma_{p}-c_{s+1,n-1}\alpha_p\gamma_{p} +
$$
$$
+\ldots+(-1)^{n-1}c_{s,1}\gamma_{p}^{n-1}+(-1)^{n-1}c_{s+1,1}\alpha_p\gamma_{p}^{n-1}=0
\eqno (14)
$$
for all $p\in P_0$\,, $s\in\{1,\ldots, n-1\}$

We consider two cases.

\textbf{Case 1.} $G'=G$. 

In this case, we choose the basis $e_1',\ldots,e_n'$ of the space $V$ and bases $v_{p,1},\ldots, v_{p,n-1}$ of the spaces $\delta_p(G)$ for $p\in P$ coinciding with the basis $e_1,\ldots,e_n$ and the basis $u_{p,1},\ldots, u_{p,n-1}$, $p\in P$, respectively.

Let $p\in P_0$, $s\in \{1,\,\ldots\,n-1\}$.
Then $\gamma_p=\alpha_p$ and relation $(14)$ takes the form
$$
d_0+d_1\alpha_p+\ldots+d_n\alpha_p^n=0,\quad \text{where}
$$
$$
d_0=c_{s,n},\quad d_1=c_{s+1,n}-c_{s,n-1},\quad 
d_2=(-1)(c_{s+1,n-1}-c_{s,n-2})
$$
$$
\cdots\cdots\cdots\cdots\cdots\cdots\cdots\cdots\cdots\cdots \eqno (15)
$$
$$
d_{n-1}=(-1)^{n-2}(c_{s+1,2}-c_{s,1}), \quad
d_n=(-1)^{n-1}c_{s+1,1}.
$$
Since the set $\{\alpha_p\,|\,p\in P_0\}$ satisfies condition $(T)$ for $n$ and $d_i\in\mathbb{Q}$, $i=0,\;\ldots,\;n$, we have  $d_0=\ldots=d_n=0$. It follows from $(15)$ that
$$
c_{s,n}=0,\; c_{s+1,1}=0,\; c_{s,r}=c_{s+1,r+1} \text{ for all } s,r\in \{1,\;\ldots,\;n-1\}. \eqno (16)
$$
Therefore, if $s<r<n$, then we have $c_{s,r}=c_{s+1,r+1}=\ldots=c_{n-(r-s),n}=0$, by $(16)$. If $s>r>1$, then $c_{s,r}=c_{s-1,r-1}=\ldots=c_{s-r+1,1}=0$; in addition $c_{1,1}=c_{2,2}=\ldots=c_{n,n}$.

Consequently, $\varphi(x)=c_{1,1}x$ for all $x\in G$, whence $\text{End }G\subseteq \mathbb{Q}$. By \cite[Lemma 2.1.8]{Fei83}, $G$ is a $nil$-group.

\textbf{Case 2.} $G'\ne G$ and the sets $\{\alpha_p\,|\,p\in P_0\}$ and $\{\gamma_p\,|\,p\in P_0\}$ are chosen so that $\alpha_p$, $\gamma_{p}$ are algebraically independent over $\mathbb{Q}$ for at least one $p\in P_0$ (such sets exist, by Lemma 14).

In this case, it follows from $(14)$ that $c_{s,i}=0$ for all $s,i\in \{1,\;\ldots,\;n-1\}$, whence $\varphi=0$. Consequently $\text{Hom }(G,G')=0$ and the groups $G$ and $G'$ are not quasi-isomorphic to each other. 

By Lemma 14, there exist $2^{\aleph_0}$ sets $\{\alpha_p\in \widehat{\mathbb{Z}}_p\,|\,p\in P_0\}$ which satisfy condition $(T)$ for $n$ and each two of which have algebraically independent $p$-components for at least one $p\in P_0$. Therefore, the class $\mathcal{MH}\,\cap\,NIL$ has $2^{\aleph_0}$ pairwise non-quasi-isomorphic groups of rank $n$ and type $t$.~$\square$

For every polynomial $f(x)\in \mathbb{Z}[x]$, we set $P(f)=\{p\in P\,|\,f(x)$ has a root in $\mathbb{Q}_p^*\}$.

\textbf{Remark 16.} It follows from the proof of \cite[Theorem 8.4]{B-P61} that if $f(x)\in \mathbb{Z}[x]$ and $\text{deg }f(x)>1$, then the set $P(f)$ is infinite. On the other hand, if $f(x)$ is irreducible in $\mathbb{Z}[x]$ and deg~$f(x)>1$, then by \cite{Fje55}, the set $P\setminus P(f)$ is also infinite.

\textbf{Theorem 17.} Let $n\in \mathbb{N}$, $n>1$, and let $t$ be a type which is not equal to $(\infty,\infty,\ldots)$. A group $G\in\mathcal{MH}\setminus \mathcal{NIL}$ of rank $n$ and type $t$ exists if and only if $t$ is an idempotent type and $P_0(t)\subseteq P(f)$ for some irreducible polynomial $f(x)\in \mathbb{Z}[x]$ of degree $n$.

\textbf{Proof.} Let $G\in\mathcal{MH}\setminus \mathcal{NIL}$, $r(G)=n>1$, $t(G)=t$. By Theorem 10, there exists a ring of field type on $G$. By \cite[Corollary 4.9]{B-P61}, the type $t$ is an idempotent type. Since $r_p(G)+\text{dim }_{\mathbb{Q}_p^*}\delta_p(G)=n$, $p\in P$, we have that 
$$
\text{dim }_{\mathbb{Q}_p^*}\delta_p(G)
=\begin{cases}
n \text{ if }r_p(G)=0\\
n-1 \text{ if }r_p(G)=1.
\end{cases}
$$
It follows from \cite[Corollary 8.5]{B-P61} that the group $G$ is quotient divisible; in addition, there exist an irreducible polynomial $f(x)\in \mathbb{Z}[x]$ of degree $n$ and a basis $e_1,\;\ldots,\;e_n$ of the space $\bar{G}$ such that
$$
\delta_p(G)=\begin{cases}
V^{(p)},\text{ if }r_p(G)=0\\
\langle e_1+\alpha_p e_2,\;\ldots,\;e_{n-1}+\alpha_pe_n\rangle_{Q_p^*},\text{ if }r_p(G)=1,
\end{cases}
$$
where $\alpha_p$ is a root of $f(x)$ contained in $\mathbb{Q}_p^*$.

Let $p\in P_0(t)$. Then $\delta_p(G)\ne V^{(p)}$. Consequently, $f(x)$ has a root $\alpha_p\in \mathbb{Q}_p^*$, i.e., $p\in P(f)$. Therefore, $P_0(t)\subseteq P(f)$.

Now let $n\in \mathbb{N}$, $n>1$,~ $t$ be an idempotent type which is not equal $(\infty,\infty,\ldots)$, and let $P_0(t)\subseteq P(f)$ for some irreducible polynomial $f(x)\in \mathbb{Z}[x]$ of degree $n$. Then for every $p\in P_0(t)$, the polynomial $f(x)$ has a root $\alpha_p\in \mathbb{Q}_p^*$.

It is easy to see that the set $\{\alpha_p\,|\,p\in P_0(t)\}$ satisfies condition $(T)$ for $n-1$. Indeed, let 
$g(x)\in \mathbb{Z}[x]$, $\text{deg }g(x)\le n-1$. Then for every $p\in P_0(t)$, the element $\alpha_p$ is not a root of $g(x)$, since $f(x)$ is a minimal polynomial for $\alpha_p$. If $h(x)$ is the largest common divisor of the polynomials $f(x)$ and $g(x))$ in the factorial ring $\mathbb{Z}[x]$, then $h(x)=c\in \mathbb{Z}$, since the polynomial $f(x)$ is irreducible. For every $p\in P_0(t)$, the relation $g(\alpha_p)\equiv 0 \text{ (mod }p)$ holds if and only if $c=h(\alpha_p)\equiv 0 \text{ (mod }p)$. Consequently, $g(\alpha_p)\not\equiv 0 \text{ (mod }p)$ for almost all $p\in P_0(t)$.

By Lemma 13, there exists a group $G\in \mathcal{MH}$ of rank $n$ and type $t$ such that for some basis $e_1,\;\ldots,\;e_n$ of the space $\bar{G}$, we have
$$
\delta_p(G)=\begin{cases}
V^{(p)} \text{ for }p\in P_{\infty}(t)\\
\langle e_1+\alpha_p e_2,\;\ldots,\;e_{n-1}+\alpha_pe_n\rangle_{\mathbb{Q}_p^*} \text{ for }p\in P_{0}(t).
\end{cases}
$$
By \cite[Corollary 8.5]{B-P61}, there exists a ring of field type on the group $G$. Therefore, $G\in \mathcal{MH}\setminus \mathcal{NIL}$.~$\square$

It follows from Theorem 17 and Remark 16 that there does not exist a homogeneous Murley group of zero type $(0;0;,\ldots)$ which admits a non-zero multiplication. We can formulate a more general assertion.

\textbf{Corollary 18.} If $t$ is a type and the set $P_{\infty}(t)$ is finite, then the class $\mathcal{MH}\setminus \mathcal{NIL}$ does not contain a group of type $t$.


\begin{thebibliography}{99}

\bibitem{AlbBV07} Albrecht U., Breaz S., Vinsonhaler C., Wickless W. Cancellation properties for quotient divisible groups // J.~Algebra. 2007. Vol.~317, no.~1. P.~424--434.

\bibitem{And67} Andrijanov V. I. Periodic Hamiltonian rings //
Mathematics of the USSR-Sbornik, 1967, Vol.~3. no.~2. P.~225--242.

\bibitem{A-P88} Andruszkiewicz R., Puczylowski E. On filial rings // Portugal. Math. 1988. Vol.~45, no.~2. P.~139--149.

\bibitem{AndS07} Andruszkiewicz R., Sobolewska M. Commutative reduced filial rings // Algebra Discrete Math. 2007. no.~7. P.~18--26.

\bibitem{AndW14} Andruszkiewicz R., Woronowicz M. On TI groups // Recent Results in Pure and Applied Math. Polasie. 2014. P.~33--41.

\bibitem{Bae57} Baer R. Meta ideals. Conference on linear algebras. June. 1956 // Publ. National Acad. Sci. nat. 1957. no.~502. P.~33--52.

\bibitem{Be48} Beaumont R. A. Rings with additive groups which is the direct sum of cyclic groups // Duke Math. J.~1948. Vol.~15, no.~2. P.~367--369.

\bibitem{B-P61} Beaumont R. A., Pierce R. S. Torsion-free rings // Illinois J. Math. 1961. Vol.~5. P.~61--98.

\bibitem{DeG57} De Groot J. Indecomposable abelian groups // Proc. Ned. Akad. Wetensch. 1957. Vol.~60. P.~137--145.

\bibitem{Er83} Ehrlich G. Filial rings // Portugal. Math. 1983-1984. Vol.~42. P.~185--194.
 
\bibitem{Fei83} Feigelstock S. Additive Groups of Rings, Pitman, Boston-London-Melbourne, 1983.

\bibitem{F-P04} Filipowicz M., Puczylowski E. R. Left filial rings // Algebra Colloq. 2004. Vol.~11. P.~335--344.

\bibitem{F-P05} Filipowicz M., Puczylowski E. R. On filial and left filial rings // Publ. Math. Debrecen 2005. Vol.~66. P.~257--267.

\bibitem{Fje55} Fjellstedt L. Bemerkungen \"uber gleichzeitige l\"osbarkeit von Kongruenzen // Arkiv Mat. 1955. Vol.~3. no.~2. P.~193--198.

\bibitem{FinR07} Fine B., Rosenberger G. Number Theory, Birkh\"auser, Boston-Basel-Berlin, 2007.

\bibitem{Fom09} Fomin A. A. Invariants for Abelian groups and dual exact sequences // J.~Algebra. 2009. Vol.~322, no.~7. P.~2544--2565.

\bibitem{Fuc15} Fuchs L. Abelian Groups, Springer Int. Publ. Switzerland, 2015.

\bibitem{Goo79} Goodearl K. R. Von Neumann Regular Rings. Pitman, London, 2002.

\bibitem{Kom10} Kompantseva E. I. Torsion-free rings // J.~Math. Sci. (New York). 2010. Vol.~171, no.~2. P.~202--211.

\bibitem{KomF19} Kompantseva E. I., Fomin A. A. Quotient divisible groups and torsion-free groups corresponding to finite Abelian groups (Russian) // Chebyshevskii sbornik. 2019. Vol.~20, no.~2. P.211--223.

\bibitem{K-T19} Kompantseva E. I., Nguyen T. Q. T. Algebraically compact abelian TI groups (Russian) // Chebyshevskii sbornik. 2019. Vol.~20, no.~1. P.~202--211.

\bibitem{KTG21} Kompantseva E. I., Nguyen T.Q.T., Gazaryan V. A. Filial rings on direct sums and direct products of torsion-free abelian groups (Russian) // Chebyshevskii sbornik. 2021. Vol.~22, no.~1. P.200--212.

\bibitem{Kru68} Kruse R. L. Rings in which all subrings are ideals // Canad. J. Math. 1968. Vol.~20. P.~862--871.

\bibitem{KMT03} Krylov P. A., Mikhalev A. V., Tuganbaev A. A. Endomorphism Rings of Abelian Groups, Springer Netherlands (Kluwer), Dordrecht–Boston–London, 2003.

\bibitem{Mur72} Murley C. E. The classification of certain classes of torsion-free abelian groups // Pacisic J. Math. 1972. Vol.~40. P.~647--665.

\bibitem{Mur74} Murley C. E. Direct products and sums of torsion-free abelian groups // Proc. Amer. Math. Soc. 1974. Vol.~38. P.~235--241.

\bibitem{Pro65} Prochazka L. A note on quasi-isomorphism of torsion-free abelian groups of finite rank // Czechoslov. Math. J. 1965. Vol.~15. P.~1--8.

\bibitem{Red52} Redei L. Vollidealringe im weiteren Sinn. I // Acta Math. Acad. 1952. Vol.~3. P.~243--268.

\bibitem{Sa88} Sands A. D. On ideals in over-rings // Publ. Math. Debrecen. 1988. Vol.~35. P.~273--279.

\bibitem{Tug02} Tuganbaev A. A. Rings Close to Regular. Springer Netherlands (Kluwer), Dordrecht–Boston–London, 2002.


\end{thebibliography}
\end{document}